\documentclass{amsart}
\usepackage{epsf,epsfig,amssymb}

\newtheorem{theorem}{Theorem}
\newtheorem{prop}[theorem]{Proposition}
\newtheorem{lemma}[theorem]{Lemma}
\newtheorem{corollary}[theorem]{Corollary}

\newcommand{\vpic}[1]{{\rule[-.7em]{0em}{2em}\:\vcenter{\epsffile{#1.eps}}\:}}
\newcommand{\Z}{{\bf Z}} 
\newcommand{\Q}{{\bf Q}} 
\newcommand{\R}{{\mathcal R}} 
\newcommand{\p}[1]{{\mathcal P}_{#1}} 
\newcommand{\V}[2]{{\mathcal V}^{#1}_{#2}}
 
\newcommand{\vpicstdheight}[1]{\kern-10pt\hbox to20pt{$\vcenter{\epsfig{file=#1.eps,width=15pt}}$}\kern 10pt}

\title[Finite Type Invariants of Pure Braids]{Free Groups and Finite Type 
Invariants \\of Pure Braids}

\author{Jacob Mostovoy}
\address{Instituto de Matem\'{a}ticas (Unidad Cuernavaca),
Universidad Nacional Aut\'{o}noma de M\'{e}xico,
A.P. 273-3 Cuernavaca, Morelos, MEXICO}
\email{jacob@matcuer.unam.mx}
\author{Simon Willerton}
\address{Institut de Recherche Math\'ematique Avanc\'ee, Universit\'e
Louis Pasteur et CNRS, 7 rue Ren\'{e} Descartes, 67084 Strasbourg Cedex,
FRANCE}
\email{willert@math.u-strasbg.fr}
\date{31 May, 1999}

\begin{document}

\begin{abstract}
In this paper  finite type invariants (also
known as Vassiliev invariants) of pure braids are considered from a 
group-theoretic point of view. New results include a construction of a 
universal invariant with integer coefficients based on the Magnus expansion 
of a free group and a calculation of numbers of independent invariants of each
type for all pure braid groups.
\end{abstract}
\maketitle

\section*{Introduction}

The attention knot theorists have paid to Vassiliev invariants of braids
can be explained in part by the fact that pure braids seem to be a good
model example of ``knotty objects'' for which all of the important questions 
about finite type invariants can be efficiently answered. 

However, the nice behaviour of Vassiliev invariants of 
pure braids is due to the very special
algebraic structure of pure braid groups. For instance,
as far as Vassiliev invariants are concerned 
pure braid groups are indistinguishable from products of free groups. 
This is of crucial importance, as the theory
of finite type invariants for free groups was developed by Fox half a
century ago under the name of ``free differential calculus'' \cite{Fox}. The
transition from the knot-theoretic language of ``overcrossings'',
``undercrossings'' and ``double points'' to group theory becomes
possible after the module generated by pure braids with $n$ double
points is identified with the $n$th power of the augmentation ideal of
the pure braid group $\p{k}$.  Recall that if $\R$ is a commutative
unital ring and $G$ is a group, the augmentation ideal
$JG\triangleleft\R G$ is defined to be the kernel of the augmentation
homomorphism, which sends all $g\in G$ to $1\in\R$.
This algebraicisation was implicit in the work of Stanford
\cite{Stanford1, Stanford2}; it appears quite explicitly in
\cite{Stanford3} and 
in the second author's thesis \cite{Simon:thesis}. 

The aim of the present paper is to give a brief exposition of
the applications of group-theoretic methods to the finite type invariants of
braids. Many of the results obtained are known. New results include 
a construction of a universal Vassiliev invariant based on the 
Magnus expansion and the calculation of the numbers of independent finite type 
invariants for all pure braid groups. The main technical tool, which is
the statement that the powers of the augmentation ideal cannot 
distinguish pure braid groups from products of free groups, also appears
to be new. Some of its corollaries, however, are well-known: one is the 
theorem of Falk and Randell \cite{FalkRandell} which describes the lower 
central series of the pure braid groups; another is the fact that the modules
of chord diagrams for pure braids are freely generated by non-decreasing 
diagrams. 

This is how the paper is organised.  In Section~1 the finite type
condition is translated into algebraic terms. 
Artin's notion of combing pure braids is considered in Section~2 .
In what is essentially the key theorem, the notion of combing is shown to 
extend in a suitable sense to singular braids.  The main consequence
of this is that many questions about finite type
invariants of pure braids can be reduced to questions about
products of free groups and these can often be answered easily.  This
philosophy is epitomised in Section~3 where the Magnus expansion of
free groups is used to give a universal finite type invariant of pure
braids which has integer coefficients.  In the fourth section the lower
central series of the pure braid group is used to characterise braids
indistinguishable form the trivial braid 
by finite type invariants of a given order.  Section~5 contains the
calculation of explicit formul\ae\ for the number of invariants at
each order.  The final section consists of some remarks concerning
relations with finite type invariants of knots.

About the notation: if $\alpha$ and $\beta$ are braids the product
$\alpha\beta$ will mean ``$\alpha$ on top of $\beta$''. 
It will sometimes be necessary to consider a linear extension
of a map between groups to their group rings; in such a situation we will use
the same letter for both maps and will not distinguish between the
two.  $\R$ will be a commutative, unital ring.


\section*{Acknowledgments}
 The first author thanks the Max-Planck-Institut f\"{u}r
Mathematik, Bonn for its kind hospitality during the stay at which this paper 
was written.  Some of this work formed a part of the second author's thesis
\cite{Simon:thesis} and for this he acknowledges the support of an EPSRC 
studentship.  Both of us thank Elmer Rees for wise counsel and 
for the proof of
Lemma~\ref{combinatoriallemma}.  We also thank Ted Stanford for helpful
conversations. 


\section{Finite type invariants of pure braids.}
This section consists of the definition of finite type invariants and
how this reduces to a purely group theoretic notion.  

In the theory of finite type invariants one considers singular pure
braids.  A singular pure braid is a pure braid whose strands are
allowed to intersect transversally at a finite (possibly
zero) number of ``double points''.  A pure braid invariant $v:\p{k}\to
\R$ can be
extended inductively to an invariant of singular pure braids by means
of the Vassiliev skein relation:
\[ v(\vpicstdheight{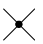}) = v(\vpicstdheight{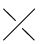}) -v(\vpicstdheight{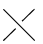}),\]
where as usual the diagrams represent pure braids which are identical
outside of some ball, inside of which they differ as shown.  An invariant
of pure braids is said to be of type $n$ if its extension vanishes on
all singular pure braids with more than $n$ double points.  An
invariant is said to a {\em finite type\/} or   {\em Vassiliev\/}
invariant if it is of type $n$ for some $n$.

For pure braids, the study of finite type invariants can
be reduced to an algebraic problem in the following manner.  A
singular braid with $k$ strands can be formally considered as an
element of the group algebra $\R\p{k}$ by the identification:
\[ \vpicstdheight{dpt} = \vpicstdheight{pos} -\vpicstdheight{neg}
 \in\R\p{k},\]
in which case the extension of an invariant 
to singular braids is precisely the same
thing as the linear extension of the invariant to the group algebra.
The investigation of finite type invariants is aided by two important
facts. 

Firstly, any pure braid can be transformed into the trivial braid by a
sequence of  
crossing changes. This implies that the augmentation ideal,  
$J\p{k}$, of the pure braid group is spanned over $\R$ by singular
braids, as
\[\begin{array}{rcl}
J\p{k} & = & \langle p-1 \ |\ p\in\p{k} \rangle\\
 & = & \langle (p_1-p_2)+(p_2-p_3)+\ldots + (p_j-1)\ | \\
 & & \qquad p_i,p_{i+1} {\rm \ differ\ by\ a\ 
crossing\ change\ (with\ } p_{j+1}=1)\rangle\\
 & = & \langle (p-q)\ |\ p,q {\rm\ differ\ by\ a\ crossing\ change}\rangle\\ 
 & = & \langle {\rm singular\ braids}\rangle.
\end{array} \] 

Secondly, any singular braid with $n$ double points can be written as a product
of $n$ singular braids, each with one double point.%
\footnote{This is where the theory differs fundamentally from that in
the case of knots.}
For example:
\[ \epsffile{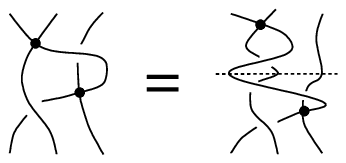} \]
This means that the $n$-th power of the augmentation ideal, $J^{n}\p{k}$, 
is spanned by singular braids with $n$ double points.

Thus, denoting by $\V{n}{k}$ the $\R$-module of type $n$ 
invariants of $\p{k}$, there is a canonical isomorphism of $\R$-modules
\[\V{n}{k}={\rm Hom}(\R\p{k}/J^{n+1}\p{k},\R).  \] 
In this sense the object of study has been reduced to something
entirely group-theoretic.

\section{Pure braids and free groups.}
In this section Artin's notion of combing braids is recalled and the
main theorem is that this extends in a suitable sense to singular
braids.  
This is used to
identify, as $\R$-modules, the powers of the augmentation ideals of pure braid
groups and those of products of free groups.

\subsection{Combed braids and combed singular braids}

Forgetting the $k$th strand of a braid gives a homomorphism  
$\p{k}\to\p{k-1}$. The kernel of this forgetful map consists of
braids which may be drawn so that the first $k-1$ strands 
are vertical and the final strand moves between them.  This can be
identified with the free group $F_{k-1}$ on $k-1$ generators.

The map $\p{k-1}\rightarrow\p{k}$ which, to a pure braid on $k-1$
strands, 
just adds a vertical, non-interacting strand on the right, is a
section of the forgetful map above.
This means that there is a split extension
\[ 1\to F_{k-1}\to\p{k}\leftrightarrows\p{k-1}\to 1,\]
or in other words that $\p{k}$ is a semi-direct product
$F_{k-1}\ltimes\p{k-1}$.  Thus, inductively there is an isomorphism
\[\p{k}\cong F_{k-1}\ltimes\ldots F_{2}\ltimes F_1.\]

\begin{figure}
\[x_{m,i}=\hspace{-2in}\vpic{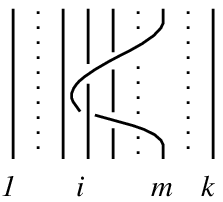}\hspace{-2in} \quad
x_{m,i}-1=\hspace{-2in}\vpic{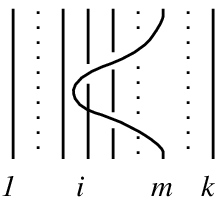}\hspace{-2in} \quad
  \sigma_i=\hspace{-2in}\vpic{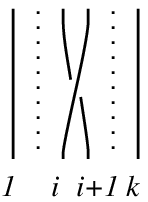} \]
\caption{Some elementary braids.}
\label{pictures}
\end{figure}

From now on $F_{m-1}$, the free group on $m-1$ generators,
will be identified 
with the free subgroup of $\p{k}$ formed
by pure braids which can be made to be totally straight apart from the
$m$th strand which is allowed to braid around the strands to the left.
Such braids   will be called $(m-1)$-free.  $F_{m-1}$
has free generators $x_{m,i}$ for $1\leqslant i \leqslant m-1$ as
illustrated 
in Figure~\ref{pictures}.  By an identical argument to that of the
last section, it can be shown that the $n$th power of the augmentation
ideal, $J^nF_{m-1}$, is spanned by $m$-free singular braids with $n$
double points.  For example, the element $x_{m,i}-1\in J^1F_{m-1}$ is
pictured as a singular braid in Figure~\ref{pictures}.

The 
semi-direct product structure means that every pure braid can be 
uniquely written as a product $\beta_{k-1}\beta_{k-2}\ldots\beta_{1}$,
where $\beta_{i}\in F_{i}$ is $i$-free. This is Artin's combing of a
pure braid \cite{Artin:TheoryOfBraids}. 

Quite similarly one can introduce the notion of a combed singular braid.
A singular braid is said to be combed if
it is written as a product $B_{k-1}B_{k-2}\ldots B_{1}$ where each
singular braid $B_{i}$ is $i$-free.  All singular pure braids can be combed
in the following sense.

\begin{theorem}\label{thm:comb}
Every singular pure braid with $n$ double points 
is equal to a linear combination of combed
singular pure braids each with $n$ double points.
\end{theorem}
This is proved in Section~2.2 below.  Theorem~\ref{thm:comb} can   be
reformulated in the following manner.  Defining the direct product of free
groups as $\Pi_{k}{:=}F_{k-1}\times\ldots
\times F_{1}$, the combing map%
\footnote{The symbol  ${\sf E}$ is supposed to evoke the image of a comb.}
${\sf E}\colon \p{k}\to \Pi_{k}$ is constructed by simply sending a
braid to its combed form, i.e.\ it is the set theoretic bijection underlying
the iterated semi-direct product structure.
The linear extension of  ${\sf E}$ gives an isomorphism of $\R$-modules
\[ {\sf E}: \R\p{k}\stackrel{\cong}{\longrightarrow}\R\Pi_{k}.\] 
Now $\R\Pi_{k}\cong \R F_{k-1}\otimes\ldots\otimes\R F_{1}$ and 
the $n$th power $J^n\Pi_k$ of the augmentation ideal of $\Pi_k$ can
be seen to be $\bigoplus_{\sum n_i =n} J^{n_{k-1}}F_{k-1}
 \otimes \dots  \otimes J^{n_{1}}F_{1}$;  so 
Theorem~\ref{thm:comb} can be formulated as follows:
\begin{theorem}\label{thm:comb2} For each $n\geqslant 0$, the combing
map induces a bijection of $n$th powers of the augmentation ideals:
${\sf E}(J^n\p{k})=J^{n}\Pi_{k}$.
\end{theorem}

Theorem~\ref{thm:comb2} is an efficient tool for obtaining information about
Vassiliev braid invariants from well-known facts about free groups.
Note, however, that the behaviour of invariants under the braid 
multiplication is not described by this theorem, as ${\sf E}$ is not a 
ring homomorphism.


\subsection{Proof of Theorem~\ref{thm:comb}} 
The proof consists of induction on the number of strands. 
The theorem is immediate for $\p{2}$. For the general case it
suffices to show that any 
singular braid on $k+1$ strands with $n$ double points can be 
written as $\sum P_{j}Q_{j}$ where $P_{j}$ is $k$-free
with $n_j$ double points and
$Q_{j}$ has  $n-n_j$
double points and has the final strand  non-interacting.

Let $B$ be a singular pure braid on   $k+1$ strands 
with $n$ double points. 
It can be written as a product $B=A_{1}s_{1}A_{2}\ldots A_{l}s_{l}$
where each $s_i$ is a singular braid with the final strand not
interacting and each $A_i$ is a $k$-free singular braid.  
It suffices to ``push the $A_i$ to the left''
in this expression for $B$, preserving the number of double points.  
Thus the theorem follows from the
following commutation lemma:
\begin{lemma}\label{lemma:pass}
If $A\in J^{m}F_{k}$ and $\sigma$ is equal to $\sigma_{i}^{\pm 1}$ where
$\sigma_{i}$ is a standard generator of the full braid group
with $1\leqslant i<k$, as in Figure~\ref{pictures}, then:\\
$\mathrm{(a)}$ there exists $A'\in J^{m}F_{k}$ such that 
$\sigma A=A'\sigma$;\\
$\mathrm{(b)}$ there exist $A'\in J^{m}F_{k}$ 
and $A''\in J^{m+1}F_{k}$ such that
$(\sigma-\sigma^{-1})A=A'(\sigma-\sigma^{-1})+A''\sigma$.
\end{lemma}
\begin{proof}[Proof of Lemma~\ref{lemma:pass}]
Part (a) is established by taking $A'=\sigma A\sigma^{-1}\in J^{m}F_{k}$ as
the powers of the augmentation ideal are invariant under all
automorphisms of $F_{k}$.

To simplify the notation, write $x_j:=x_{k,j}$ for the generators of
the free group $F_{k}$.  By commuting generators and double points
one at a time, it is sufficient to verify (b) for the 
generators $\{ x_j\}$, their inverses, and the set of elementary
singular braids   $\{ x_j-1\}\in JF_{k}$.
It is clear that if $j\neq i,i+1$ the generator $x_{j}$ commutes with 
$\sigma$.
Hence, it remains to verify the lemma for $x^{\pm 1}_{i}$ and 
$x^{\pm 1}_{i+1}$. We will do the calculations only for $x_{i}$ as the case of 
$x_{i}^{-1}$ and $x^{\pm 1}_{i+1}$ can be treated in the same manner.

The following relations are easily checked:
$$\sigma x_{i}=x_{i+1}\sigma,\qquad
\sigma^{-1} x_{i}=x_{i} x_{i+1} x_{i}^{-1}\sigma^{-1}.$$
Thus,
\[\begin{array}{ccl}
(\sigma-\sigma^{-1})x_{i}& = & x_{i+1}\sigma-x_{i}^{-1}
x_{i+1}x_{i}\sigma^{-1}\\
  & = &  x_{i}^{-1}x_{i+1}x_{i}(\sigma-\sigma^{-1})+
(x_{i+1}-x_{i}^{-1}x_{i+1}x_{i})\sigma\\
  & = & Y'(\sigma-\sigma^{-1})+Y''\sigma,
\end{array}\]
where $Y'\in F_{k}$ and
$$Y''  =  x_{i+1}-x_{i}x_{i+1}x_{i}^{-1} 
       =  x_{i}(x_{i+1}-1)(1-x_{i}^{-1})-(x_{i}-1)(x_{i+1}-1)
       \in J^{2}F_{k}.$$

Finally,  $(\sigma-\sigma^{-1})(x_{i}-1)=(Y'-1)(\sigma-\sigma^{-1})+Y''\sigma$
where $Y'$ and $Y''$ are as above, so 
$(Y'-1)\in JF_{k}$ and $Y''\in J^{2}F_{k}$. 
\end{proof}


\section{The Magnus expansion as a universal Vassiliev invariant.}

An explicit universal Vassiliev invariant for pure braids 
is provided by Chen's expansion \cite{Kohno}, 
also known in this context as the Kontsevich
integral. The Kontsevich integral is multiplicative and its value on the 
generators of the full braid group is easily calculated. Nevertheless,
finding the Kontsevich integral of an arbitrary braid is a non-trivial problem;
one difficulty resides in expressing the answer in terms of some fixed basis
for the vector space of chord diagrams.

The coefficients of the Kontsevich integral are rational numbers
rather than integers;
but universal invariants with integral coefficients are also known to exist.
A construction which gives such invariants was 
described by Hutchings in \cite{Hutchings}. Hutchings' method 
yields all universal invariants; however, it is not very practical. 
In this section a universal braid invariant with integer coefficients 
which can be easily calculated, is constructed using the
Magnus expansion.   Other expansions of the free 
group (see  \cite{Lin}) also give universal pure braid
invariants; the resulting universal invariant has integer coefficients
if the corresponding expansion of the free group does.


\subsection{The Magnus expansion for a free group.} 
The statements below are essentially contained in \cite{Fox}.

Let $G$ be a group and set  $A_{n}(G)=J^{n}G/J^{n+1}G$. Define the $\R$-module 
\[ A(G):=\widehat\bigoplus A_{n}(G) \]  
to be the completion of the direct
sum with respect to the grading. The module $A(G)$ is, in fact, an 
algebra, with multiplication induced by that on $G$. Sometimes  $A^{\R}(G)$ 
will be written to emphasise the ground ring. In the case that $G$ is a
pure braid group, $A(\p{k})$ is often called the {\em algebra of chord 
diagrams.}

For $G$ a free group this algebra has a particularly simple structure.
Namely, let $\R[[X_{1},\dots,X_{i}]]$ be the algebra of formal power series 
in $i$ non-commuting variables. There is an isomorphism of algebras 
\[ \R[[X_{1},\dots,X_{i}]]\stackrel{\cong}{\longrightarrow} A(F_{i}) \]
which sends $X_{j}$ to $x_{j}-1$, where $x_{j}$ are the generators of
$F_{i}$. In particular, the abelian group $A^{\Z}_{n}(F_{i})$ is 
torsion-free and its rank is equal to the number of different monomials
of degree $n$ in $i$ non-commuting variables, i.e.\ to $i^{n}$.

The {\em Magnus expansion\/} is the algebra homomorphism
\[ M\colon\R F_{i} \to \R[[X_{1},\dots,X_{i}]] \]
defined by sending the generator $x_{j}$ to $1+X_{j}$ for 
$1\leqslant j \leqslant i$. This definition implies for example, that
$M(x_{j}^{-1})=1-X_{j}+X_{j}^{2}-X_{j}^{3}+\dots$.  
A fundamental property is the following. 
\begin{prop}\label{prop:magnusfree}
If $x\in J^{n}F_{i}$ then $M(x)$ has no terms of degree less than $n$
and the $n$th degree term is the image of $x$ under the natural projection
\[ J^{n}F_{i}\to A_{n}(F_{i})\cong(\R[[X_{1},\dots,X_{i}]])_{n}.\]
\end{prop}


\subsection{The Magnus expansion for pure braid groups} 
The algebra $A(\p{k})$ is often called 
the algebra of chord diagrams, and by
Theorem~\ref{thm:comb2}  there is an isomorphism of $\R$-modules
\[\widetilde{\sf E}: {A}(\p{k})\stackrel{\cong}{\longrightarrow} {A}(\Pi_{k})
\cong {A}(F_{k-1})\otimes{A}(F_{k-2})\otimes\dots\otimes{A}(F_{1}).\]
This implies that ${A}^\Z_n(\p{k})$ is torsion free, and also 
provides ${A}_n(\p{k})$ with a basis. Indeed, if
$A(F_{i})$ is the non-commutative power series ring in the degree one
indeterminates  $ \{X_{i+1,j}\}$, where $j\leqslant i$, then a basis for
${A}_n(\p{k})$ is given by monomials of the form 
$m_{k-1}\otimes\ldots\otimes m_{1}$ where each $m_{i}\in
A_{n_i}(F_{i})$ is a degree $n_i$ monomial and $\sum n_i=n$.
These can be represented graphically as the so-called {\em
non-decreasing\/}  or {\em sorted\/} diagrams: $k$ vertical strands
are drawn and  each $X_{p,q}$ is represented as a 
horizontal chord between the $p$th and $q$th strand.  For example, a
basis element of ${A}_5(\p{3})$ is drawn:
$$X_{3,1}X_{3,2}X_{3,1}\otimes X_{2,1}X_{2,1}
\longleftrightarrow \vpicstdheight{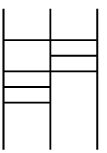}.$$

The Magnus expansion of pure braids can be defined as follows.  Let $\beta$
be a combed braid: $\beta=\beta_{k-1}\beta_{k-2}\ldots\beta_{1}$,
with $\beta_{i}\in F_{i}$. Define the Magnus expansion 
$M\colon \R\p{k}\to A(\p{k})$ by
\[ M(\beta):=\widetilde{\sf E}^{-1}\bigl(M\left(\beta_{k-1}\right)
\otimes\ldots\otimes M\left(\beta_{1}\right)\bigr).\]  

The Magnus expansion of a pure braid can be easily computed in the basis of 
non-decreasing diagrams, as the problem reduces to calculations of Magnus 
expansions in free groups. This is best illustrated with an example:
\begin{align*}
M\Biggl(\vpicstdheight{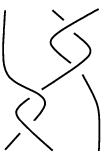}\Biggr)
&=\Biggl( 1+\vpicstdheight{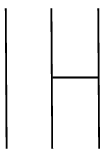} \Biggr)
  \Biggl( 1-\vpicstdheight{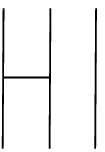}+ \vpicstdheight{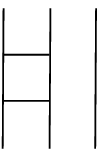}
              -\vpicstdheight{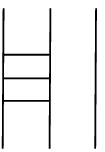}+\ldots  \Biggr)\\
&=1+\vpicstdheight{diag1}-\vpicstdheight{diag2}-
\vpicstdheight{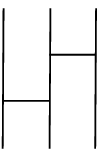}+\vpicstdheight{diag4}+
\vpicstdheight{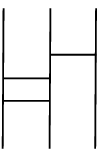}-\vpicstdheight{diag6}+
\ldots.
\end{align*}

This illustrates two very nice properties of this invariant: the
first being that its value is naturally expressed in terms of
a basis of chord diagrams and the second is that it has integer
coefficients.  Unfortunately there are prices to be paid for these
niceties: the first being that to calculate the invariant the braid
has to be initially combed and the second is that it is not an algebra
map, unlike the Kontsevich integral.  However, the latter is not
surprising as any universal Vassiliev invariant with integer
coefficient lacks this property, see \cite{Hutchings}.

Note though that this map is multiplicative with respect to the
external products $\p{k}\times \p{l}\to \p{k+l}$ and $A(\p{k})\otimes
A(\p{l}) \to A(\p{k+l})$, defined in both cases by placing a
$k$-strand object to the left of an $l$-strand object.
The reader is 
left to work out the details.

\subsection{The Magnus expansion as a universal Vassiliev invariant.}
In the context of pure braids a {\em universal Vassiliev invariant\/} with 
coefficients in $\R$ is an $\R$-linear map  $U\colon\R\p{k}\to A^{\R}(\p{k})$ 
such that for any $\R$-valued 
finite type invariant $v$ there is a unique $\R$-linear map 
$\hat{v}\colon A(\p{k})\to\R$ such that $\hat{v}\circ U=v$.

\begin{theorem}\label{thm:magnusbraids}
The Magnus expansion is a universal Vassiliev invariant.
\end{theorem}
\begin{proof}
In view of Theorem~\ref{thm:comb2} it is clear that 
Proposition~\ref{prop:magnusfree} implies 
an analogous statement for the Magnus expansion of braids. Namely,
for any braid $\beta$ with at least $n$ double points $M(\beta)$ has no
terms of  degree less than $n$ and the $n$th degree term is the image
of $\beta$ under the natural projection $J^{n}\p{k}\to A_{n}(\p{k})$.

The proof proceeds by induction on the type of an invariant $v$.  
Invariants of type zero are just constants and the existence of a 
homomorphism $\hat{v}: A(\p{k})\to \R$ such that $\hat{v}\circ M = v$ for 
them is 
trivial. Suppose that homomorphisms of $A(\p{k})$ to $\R$ which correspond to 
invariants of type less than $n$ have been found.
Let $v$ be of type $n$. As $v$ vanishes on $J^{n+1}\p{k}$, it defines
a homomorphism $\hat{v}_{1}:A_{n}(\p{k})\to\R$ and can be trivially extended
to a homomorphism 
$\hat{v}_{1}:A(\p{k})\to\R$. Then $\hat{v}_{1}\circ M - v$ is an
invariant of type $n-1$ 
and, hence by the inductive hypothesis there exists a $\hat{v}_{2}$
such that $\hat{v}_{1}\circ M - v=\hat{v}_{2}\circ M.$ Now, setting 
$\hat{v}=\hat{v}_{1}+\hat{v}_{2}$ we obtain a homomorphism
$\hat{v}:A(\p{k})\to\R$ which represents $v$. The uniqueness is 
straightforward.
\end{proof}

As the Magnus expansion  over $\Z$ for free groups is known to
be injective, there is an immediate corollary:
\begin{corollary}
Integer-valued finite type invariants separate braids.
\end{corollary}
\noindent This was also obtained in \cite{Kohno:YB} and \cite{Bar:VasQIB}.

The very existence of a universal Vassiliev invariant with coefficients in 
$\R$ has an important consequence.
The $\R$-module of all Vassiliev invariants is 
filtered by the type of the invariants:
$\V{1}{k}\subset \V{2}{k}\subset \dots$. From Section~1 there is the
natural identification $\V{n}{k}/\V{n-1}{k}\cong {\rm Hom}(A_{n},\R)$,
of the $\R$-module of type $n$ invariants modulo the type $n-1$ invariants 
with the dual of the $\R$-module $A_{n}$.  The existence of a universal 
Vassiliev invariant gives a splitting of the
filtration, i.e.\  an isomorphism from the module of Vassiliev
invariants to its associated graded module $\bigoplus {\rm Hom}(A_{n},\R)$.
Thus the module of Vassiliev invariants is dual to the algebra $A(\p{k})$.

As the modules $A_{n}^{\Z}(\p{k})$ are torsion-free, one corollary is
that all $\Z/p$-valued finite type invariants
of pure braids are just 
mod $p$ reductions of integer-valued invariants of the same type.


\section{Lower central series and {\em n}-triviality.}
An $n$-trivial pure braid is one which cannot be distinguished from
the trivial braid by invariants of type less than $n$.
In this section an easy characterisation of $n$-trivial pure braids is
given.   The $n$-triviality of pure braids was considered by Stanford
\cite{Stanford1,Stanford2,Stanford3}. 
Here we follow \cite{Stanford3} and \cite{Simon:thesis}.

The $n$th group $\gamma_{n}G$ of the {\em lower central series\/}
of a group $G$ is
defined inductively by setting $\gamma_{1}G:=G$
and, for $n > 1$, is 
defined as the subgroup generated by commutators of elements of
$\gamma_{n-1}G$ with elements of $G$:
 \[\gamma_{n}G:=[\gamma_{n-1}G, G].\]
These subgroups form a descending filtration of $G$
\[ G=\gamma_1 G\triangleright \gamma_2 G\triangleright \dots \triangleright
\gamma_n G \triangleright  \cdots.\]
Each $\gamma_{n+1}G$ is normal in  $\gamma_{n}G$ and the quotients
$\gamma_{n}G/\gamma_{n+1}G$ are abelian. 

The {\em dimension subgroups\/} are defined for $n\geqslant 1$ by
\[ \Delta_{n}G=(1+J^{n}G)\cap G,\]
note that this depends on the choice of the ring $\R$.  If $G$ is a
pure braid group then the $n$th dimension subgroup, $\Delta_{n}\p{k}$,
consists precisely of the $n$-trivial braids.

It is not hard to check that the groups of the lower central series are
contained in the dimension subgroups, 
$\gamma_{n}G\triangleleft\Delta_{n}G$. If
equality holds then then the group $G$ is said 
to have the {\em dimension subgroup property\/}.%
\footnote{It was an open question for 
many years if there is a group which does {\em not} have this property. 
The example of such a group was given by Rips in \cite{Rips}.}

Free groups are known to have the dimension subgroup property
\cite{Hartley} and a theorem of 
Sandling \cite{Sandling} says that a semi-direct product of groups that
have the dimension subgroup property has this property itself; these
facts imply that the pure braid groups have the dimension subgroup property.
An immediate consequence is the following, which was also shown by
Kohno \cite{Kohno:Contemp}.
\begin{theorem}
The $n$-trivial braids are precisely those in the $n$th subgroup,
$\gamma_n\p{k}$, of the lower central series.
\end{theorem}

Another consequence of the dimension subgroup property of $\p{k}$ and
Theorem~\ref{thm:comb2} is that the combing map respects the lower
central series:
\[ {\sf E}(\gamma_{n}\p{k})=\gamma_{n}\Pi_{k}. \]
This was proved by Falk and Randell in \cite{FalkRandell}. Note that
$\gamma_{n}\Pi_{k}$ is a direct product of $\gamma_{n}F_{i}$ for all
$i<k$ as the lower central series of a direct product is readily seen
to be the product of the lower central series of the factors.


\section{Counting the numbers of invariants.}

The purpose of this section is to obtain the number of linearly
independent invariants of
each type for each pure braid group. Throughout the section it is assumed that
$\R=\Q$.

As the Vassiliev invariants form the dual of the algebra $A(\p{k})$,
calculating the numbers of linearly independent invariants amounts
to finding the dimensions of each $A_{n}(\p{k})$. Such a computation was 
done in \cite{100} where the dimension of  $A_{n}(\p{k})$ was shown to be
a Stirling number of the second kind.

However, among Vassiliev invariants are some that are equal to sums of 
products of invariants of lower order. It seems sensible to exclude these from 
consideration --- more precisely, to factor them out --- and count only 
indecomposable invariants; this is done in Section~5.1. 
A further reduction is possible if we consider 
only those invariants that are not induced from pure braid groups on fewer 
strands; this is done in Section~5.2.


\subsection{Finding the numbers of indecomposable invariants.}

The multiplication in $\Q$ induces a multiplication on $\Q$-valued
finite type invariants and gives the graded module 
$\bigoplus \V{n}{k}/\V{n-1}{k}={\rm Hom}(A(\p{k}),\Q)$ the structure of a
commutative algebra.   To find the number of invariants of each degree
which do not come from those of lower degree it is necessary to find
the dimension
of the vector space of indecomposable elements of this graded algebra
in degree $n$.  Let $\varphi^k_n$ be this dimension: a closed
expression for this number is given by the next theorem.

\begin{theorem}
The dimension, $\varphi^k_n$, of the space
of indecomposable type $n$ pure braid invariants modulo type $n-1$
invariants is given by
\[ \varphi^k_n=\frac{1}{n} \sum_{m\mid n}\mu \left(n/m\right) 
 \sum_{i=1}^{k-1}i^m \]
where $\mu$ is the M\"obius function of number theory.
\label{dimindecthm}
\end{theorem}
\noindent Some values of $\varphi^k_n$ are tabulated in
Table~\ref{varphitable}. 

\begin{proof}
By duality, the product on ${\rm Hom}(A(\p{k}),\R)$ gives rise to a coproduct
on $A(\p{k})$. 
This is also the coproduct induced by the natural coproduct on
$\Q\p{k}$, namely $\beta\to\beta\otimes\beta$ for each pure braid $\beta$.

From a result of Milnor and Moore \cite{MiMoo}, the space of
indecomposable elements of a commutative algebra is dual to the space
of primitive elements in the dual coalgebra, so $\varphi^k_n$ is equal
to the dimension of the subspace of primitive elements in $A_{n}(\p{k})$.

Now, denote by $G\{ n\}$ the abelian group $\gamma_{n} G/\gamma_{n+1} G$. The
graded group $\bigoplus G\{ n\}$ is a Lie algebra with the Lie bracket 
induced by the group commutator. A theorem of Quillen \cite{Quillen}
says that the algebra
$\bigoplus A_n^{\Q}(G)$ is the universal enveloping algebra of
of the Lie algebra $\bigoplus G\{ n\}\otimes\Q$. However, the
space of primitive elements of a universal enveloping algebra is 
naturally isomorphic to the original Lie algebra. This means that 
$\varphi^k_n$ is equal to the rank of the abelian group $\p{k}\{ n\}$.

The rank of $\p{k}\{ n\}$ can be found from Kohno's work \cite{Kohno}
or expressed via the ranks of $F_{i}\{ n\}$ using a result of Falk and Randell
\cite{FalkRandell}.
Alternatively, notice that if $A_{n}(\Pi_{k})$ is considered as a coalgebra 
with the coproduct induced by the natural coproduct on $\R\Pi_{k}$, the 
combing identification $\widetilde{\sf E}:A(\p{k})\to A(\Pi_{k})$ is
an isomorphism of coalgebras. Hence, the primitive elements in $A(\Pi_{k})$
and $A(\p{k})$ are in one-to-one correspondence. Now, according to Quillen's
theorem, $\varphi^k_n$ is equal to the rank of the abelian group 
$\Pi_{k}\{ n\}$.

Observe that for any groups $G_{1}$ and $G_{2}$ the abelian group
$(G_{1}\times G_{2})\{ n\}$ is isomorphic to $G_{1}\{ n\}\times G_{2}\{ n\}$.
Hence
\[{\rm rank}\ \Pi_{k}\{ n\}=\sum\nolimits_{i=1}^{k-1}{\rm rank}\ F_{i}\{ n\}.\]
The theorem then follows from Witt's formula for the ranks of the quotient
groups $F_{i}\{ n\}$, see \cite{MKS}:
\[{\rm rank}\ F_{i}\{ n\}=\frac{1}{n} \sum_{m\mid n}\mu \left(n/m\right) i^m.\]

\end{proof}

\begin{table}
\begin{center}\tiny
\begin{tabular}{|l||r|r|r|r|r|r|r|r|r|}
\hline
$k\backslash n$ &1& 2& 3& 4& 5&6& 7& 8& 9
 \\ \hline
2   &1& 1& 1& 1& 1& 1& 1& 1& 1\\
3     &3& 7& 15& 31& 63& 127& 255& 511& 1023\\
4           &6& 25& 90& 301& 966& 3025& 9330& 28501& 86526\\
5     &10& 65& 350& 1701& 7770& 34105& 145750& 611501& 2532530\\
6  &15& 140& 1050& 6951& 42525& 246730& 1379400& 7508501& 40075035\\
7  &21& 266& 2646& 22827& 179487& 1323652& 9321312& 63436373&
420693273\\
8   &28& 462& 5880& 63987& 627396& 5715424& 49329280& 408741333&
     3281882604\\
\hline
\end{tabular}
\end{center}
\caption{The dimensions of the spaces of type $n$
invariants modulo type $(n-1)$ invariants of the pure braid groups,
$\p{k}$ (included for comparison).}  
  \end{table}

\begin{table}                      
\begin{center}
\footnotesize
\begin{tabular}{|l||r|r|r|r|r|r|r|r|r|r|}
\hline
$k\backslash n$ &1& 2& 3& 4& 5&6& 7& 8& 9& 10
 \\ \hline
  
2&1& 0& 0& 0& 0& 0& 0& 0& 0& 0  \\
3&3& 1& 2& 3& 6& 9& 18& 30& 56& 99  \\
4&6& 4& 10& 21& 54& 125& 330& 840& 2240&   
  5979  \\ 
5 &10& 10& 30& 81& 258& 795& 2670& 9000&
  31360& 110733 \\
6 &15& 20& 70& 231& 882& 3375& 13830& 57750&248360& 1086981\\
7  &21& 35& 140& 546& 2436& 11110& 53820& 267540& 1368080&
   7132818\\
8 &28& 56& 252& 1134& 5796& 30654& 171468& 987840& 5851776
 &  35378658\\
\hline
\end{tabular}
\end{center}
\caption{The dimensions, $\varphi^k_n$ of the spaces of indecomposable
 type $n$ invariants modulo type $(n-1)$ invariants of the pure braid
 groups $\p{k}$ --- see Theorem~\ref{dimindecthm}.} 
\label{varphitable}
  \end{table}


\subsection{Reducing by invariants from lower pure braid groups.} 
There are $\binom{k}{l}$ maps from the $k$-strand pure braid
group to the $l$-strand pure braid group obtained by picking $l$ strands and
``forgetting'' the rest, thus each invariant of  $l$-strand braids
induces $\binom{k}{l}$ invariants of $\p{k}$.  For instance, all type one
invariants are linear combinations of winding numbers, and so are induced from
$\p{2}$.

To see how
many genuinely new invariants come from the $\p{k}$, one can
calculate the dimension, $\psi^k_n$, of the space of type $n$ indecomposable 
invariants
of $\p{k}$ modulo the type $n$ invariants induced from lower braid groups.
 Then these dimensions satisfy 
  $$\varphi^l_n = \sum^{l}_{k=1} \binom{l}{k}\psi^k_n.$$
Define $\operatorname{sur}(m,k)$ for $m>0$,
to be the number of surjections from an
$m$ element set to a $k$ element set, with the convention that
$\operatorname{sur}(m,0)=0$. 

\begin{theorem} The reduced dimensions, $\psi^k_n$, of type $n$
invariants of $\p{k}$  are given by 
$$\psi^k_n = \frac{1}{n}\sum_{m\mid n}  \mu \left( n/m
                    \right) \mathrm{sur}(m,k-1).$$
\label{redindecdimthm}
\end{theorem}
 \begin{table}
\begin{center}
  \begin{tabular}{|l||r|r|r|r|r|r|r|r|r|r|}
  \hline
   $k\backslash n$ &1& 2& 3& 4& 5&6& 7& 8& 9& 10
    \\ \hline
   2&1& 0& 0& 0& 0& 0& 0& 0& 0& 0
    \\
  3&0& 1& 2& 3& 6& 9& 18& 30& 56&
   99 \\
 4&0& 0& 2& 9& 30& 89& 258& 720&
  2016& 5583 \\
  5&0& 0& 0& 6& 48& 260& 1200& 5100&
  20720& 81828 \\
  6&0& 0& 0& 0& 24& 300& 2400& 15750&
  92680& 510288 \\
  7&0& 0& 0& 0& 0& 120& 2160& 23940& 211680&
  1643544 \\
  8&0& 0& 0& 0& 0& 0& 720& 17640&   
  258720& 2963520 \\
  \hline
  \end{tabular}
\end{center}  
\caption{The dimensions, $\psi^k_n$,
 of the spaces of reduced type $n$ indecomposable
invariants modulo type $(n-1)$ invariants of the pure braid groups
 $\p{k}$ --- see Theorem~\ref{redindecdimthm}.} 
\label{IndecompTable}
\end{table}
\noindent The key point is the following combinatorial identity ---
the proof given here is due to Elmer Rees.
\begin{lemma} If $m>0$ then
$\displaystyle\sum\limits_{i=1}^{l-1} i^m = \sum\limits_{j=1}^l \binom{l}{j}
\operatorname{sur}(m,j-1)$.
\label{combinatoriallemma}
\end{lemma}
\begin{proof}[Proof of Lemma~\ref{combinatoriallemma}]
The left hand side can be seen as the number of
maps from an $m$ element set to the following set, such that the image
is `vertical'.
\[ \epsfig{file=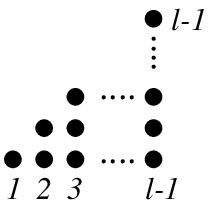} \]  
The proof proceeds by counting these maps in a different way.  So,
\begin{align*}
  \sum_{i=1}^{l-1} i^m&= \sum_{r=1}^{l-1} (\text{number of maps of image
size $r$})
    =\sum_{r=1}^{l-1} \sum_{k=1}^{l-1} \binom{k}{r}
                                 \operatorname{sur}(m,r) \\
   &=\sum_{r=1}^{l-1} \binom{l}{r+1}  \operatorname{sur}(m,r)
        =\sum_{j=1}^{l} \binom{l}{j} \operatorname{sur}(m,j-1).
\end{align*}
The third equality follows from a simple identity and the final
equality comes from relabelling and the convention that
$\operatorname{sur}(m,0)=0$. 
\end{proof}
\begin{proof} [Proof of Theorem~\ref{redindecdimthm}]
\begin{align*}
  \sum_{k=1}^{l} \binom{l}{k}\psi^k_n&=\varphi^l_n 
    =  \frac{1}{n} \sum_{m\mid n} \mu \left(
                         n/m \right) \sum_{i=1}^{l-1} i^m
    =  \frac{1}{n} \sum_{m\mid n} \mu \left( n/m \right)
                  \sum_{j=1}^l \binom{l}{j}\operatorname{sur}(m,j-1)\\
    &=   \sum_{j=1}^l \binom{l}{j} \left[  \frac{1}{n} \sum_{m\mid n}
               \mu \left( n/m \right) \operatorname{sur}(m,j-1)
                 \right].
\end{align*}
The theorem then follows as the matrix $\left( \binom{i}{j} \right)
_{1\leqslant i,j \leqslant l}$ is invertible --- it has inverse  
$\left( \left(-1\right)^{i-j} 
\binom{j}{i} \right) _{1\leqslant i,j \leqslant l}$.
\end{proof}

These dimensions are tabulated in Table~\ref{IndecompTable}.  Note
that, by \cite{Bar:VasHSL}, the entries along the leading diagonal of the 
table correspond to Milnor invariants. The fact that all invariants of type
$< k-1$ come from braids on fewer numbers of strands is not surprising
if one thinks in terms of chord diagrams:  a connected diagram on $k$ 
vertical strands has at least $k-1$ horizontal chords.


\section{Final Remarks.}
The machinery presented in this paper is very specific to pure braids.
Nevertheless, connections with knot theory do exist, some are outlined
below.


\subsection{Knots via braid closures}
The theory of finite type invariants for knots can be translated into 
the group-theoretic language by means of braid closures. A convenient
closure for this purpose is the ``short-circuit'' closure of \cite{MoSt}. 
It provides a map from the pure braid group on the infinite number of strands 
$\p{\infty}$ to the set of isotopy classes of oriented knots; this map 
can be interpreted as a two-sided 
quotient of $\p{\infty}$ by the action of two explicitly identified subgroups.
Finite type invariants of knots pull back via this map to finite type 
invariants of braids; so the problem of studying Vassiliev
knot invariants can be interpreted as the problem of studying the behaviour
of Vassiliev braid invariants under a certain two-sided action on $\p{\infty}$.

It can be shown that the short-circuit map sends the filtration of $\p{\infty}$
by the lower central series to the filtration by $n$-trivial
knots. Thus one is lead to study the interaction of commutators with 
two-sided actions of subgroups of $\p{\infty}$. Problems of this kind, however,
seem to have received little attention in group theory.


\subsection{Characteristic classes of knots}
No construction of a universal integer-valued knot invariant is known
at the moment. One may expect the relationship between such an invariant  
and the Kontsevich integral for knots to be similar in some way to that
between the Magnus expansion and the Kontsevich integral for pure braids.
The latter bears some resemblance to the relationship between the total Chern 
class and the Chern character of a vector bundle. 
Therefore, one may try to treat the Kontsevich integral of a knot as a
Chern 
character and ask if the corresponding formally defined total Chern class 
is integer-valued.

This question is considered in \cite{Simon:chern}. 
The answer turns out to be rather surprising: the total 
Chern class of a knot is 
integral, but only ``on the level of Lie algebras'': on the level
of chord diagrams this fails, the total Chern class of a trefoil being a
counterexample.


\subsection{Magnus expansion and the Gusarov-Polyak-Viro invariant}
A Magnus-type expansion has shown up in knot theory in the form of
the Gusarov-Polyak-Viro universal invariant of virtual knots
\cite{GPV}. If a knot is thought of as being ``generated by its crossings'',
then this invariant is, essentially, the ``Magnus expansion in crossings''.
For example, for a knot with only positive crossings, such as the
right-trefoil, the value of the Gusarov-Polyak-Viro  invariant 
is the formal sum of all subdiagrams of the knot
diagram. Similarly, the Magnus expansion of a word in a free group 
that contains only positive powers of generators is just the 
sum of all its subwords.

\bibliographystyle{amsplain}
\bibliography{magnus} 
\end{document}